\documentclass[12pt]{article}
\usepackage{amsfonts}
\usepackage{amsmath}
\usepackage{amsthm}
\usepackage{graphicx}
\numberwithin{equation}{section}

\newtheorem{theorem}{Theorem}[section]
\begin{document}

\title{Asymptotic zero distribution of Jacobi-Pi\~neiro and multiple Laguerre polynomials}
\author{Thorsten Neuschel, Walter Van Assche \\
Katholieke Universiteit Leuven, Belgium}
\date{\today}
\maketitle

\begin{abstract}
We give the asymptotic distribution of the zeros of Jacobi-Pi\~neiro polynomials and multiple Laguerre polynomials of the first kind.
We use the nearest neighbor recurrence relations for these polynomials and a recent result on the ratio asymptotics of multiple
orthogonal polynomials. We show how these asymptotic zero distributions are related to the Fuss-Catalan distribution.
\end{abstract}

\section{Introduction and main results}

In this paper we obtain the asymptotic distribution of the zeros of two families of multiple orthogonal polynomials: 
the Jacobi-Pi\~neiro polynomials and the multiple Laguerre polynomials of the first kind \cite[Ch. 23]{Ismail}, \cite{AptBraWVA}, \cite{WVAEC}.
These are two families of multiple orthogonal polynomials for which explicit formulas are known and which are useful for a number of applications.
For instance, the zeros of Jacobi-Pi\~neiro polynomials (and Wronskian-type determinants of Jacobi-Pi\~neiro polynomials) form the unique solution of
certain Bethe Ansatz equations \cite{MukhinVarchenko} and multiple orthogonal polynomials are also useful for investigating  determinantal
point processes \cite{Kuijl}. Recently the Jacobi-Pi\~neiro ensemble and the multiple Laguerre ensemble were introduced for random matrix minor
processes related to percolation theory \cite{AdlvMoerWang} which are based on the Jacobi-Pi\~neiro and multiple Laguerre polynomials
of the first kind. 

Let $\vec{n}=(n_1,n_2,\ldots,n_r) \in \mathbb{N}^r$ be a multi-index of size $|\vec{n}| = n_1+n_2+\cdots+n_r$.
The Jacobi-Pi\~neiro polynomials $P_{\vec{n}}$, with parameters $\vec{\alpha}=(\alpha_1,\ldots,\alpha_r)$ and $\beta$, are type II
multiple orthogonal polynomials on $[0,1]$ for $r$ Jacobi weights, i.e., $P_{\vec{n}}$ is a monic polynomial of degree $|\vec{n}|$ satisfying
\[  \int_0^1 P_{\vec{n}}(x) x^k x^{\alpha_j} (1-x)^\beta \, dx = 0, \qquad
    k=0,1,\ldots,n_j-1, \]
for $j=1,2,\ldots,r$, where $\beta > -1$ and $\alpha_j > -1$ for $1 \leq j \leq r$. They were introduced by Pi\~neiro for $\beta=0$ \cite{Pineiro}.
A multi-index $\vec{n}$ is normal if the monic multiple orthogonal polynomial $P_{\vec{n}}$ of degree $|\vec{n}|$ exists and is unique.
All multi-indices for Jacobi-Pi\~neiro polynomials are normal when $\alpha_i - \alpha_j \notin
\mathbb{Z}$ because then the measures form an AT-system \cite[\S 23.1.2]{Ismail}. The polynomials are given by the Rodrigues formula
\begin{multline} \label{JP-Rod}
 (-1)^{|\vec{n}|} \prod_{j=1}^r (|\vec{n}|+\alpha_j+\beta+1)_{n_j} \ (1-x)^\beta P_{\vec{n}}(x) \\
   = \prod_{j=1}^r \left( x^{-\alpha_j} \frac{d^{n_j}}{dx^{n_j}} x^{n_j+\alpha_j} \right) (1-x)^{|\vec{n}|+\beta},
\end{multline}
where the product of the differential operators can be taken in any order, since these operators are commuting \cite[\S 23.3.2]{Ismail}.
Multiple Laguerre polynomials of the first kind $L_{\vec{n}}$ are given by 
the Rodrigues formula 
\begin{equation}  \label{ML-Rod}
  (-1)^{|\vec{n}|} e^{-x} L_{\vec{n}}(x) = \prod_{j=1}^r \left( x^{-\alpha_j} \frac{d^{n_j}}{dx^{n_j}} x^{n_j+\alpha_j} \right) e^{-x}, 
\end{equation}
where the product of the differential operators can be taken in any order \cite[\S 23.4.1]{Ismail}.
If the parameters $\vec{\alpha} = (\alpha_1,\alpha_2,\ldots,\alpha_r)$ are such that
$\alpha_i > -1$ for every $i$ and $\alpha_i - \alpha_j \notin \mathbb{Z}$ $(1 \leq i,j \leq r)$, 
then all multi-indices are normal and the polynomials satisfy the following orthogonality properties 
\[  \int_0^\infty L_{\vec{n}}(x) x^k x^{\alpha_j} e^{-x}\, dx = 0, \qquad
     k=0,1,\ldots,n_j-1, \]
for $j=1,2,\ldots,r$.
An explicit expression is given by
\begin{multline}  \label{ML-explicit}
  L_{\vec{n}}(x) = \sum_{k_1=0}^{n_1} \cdots \sum_{k_r=0}^{n_r} (-1)^{|\vec{k}|} \frac{n_1!}{(n_1-k_1)!}\cdots \frac{n_r!}{(n_r-k_r)!} \\
   \times  \binom{n_r+\alpha_r}{k_r} \binom{n_r+n_{r-1}+\alpha_{r-1}-k_r}{k_{r-1}} \cdots  
 \binom{|\vec{n}|-|\vec{k}|+k_1+\alpha_1}{k_1}     x^{|\vec{n}|-|\vec{k}|}.
\end{multline}

We will obtain the asymptotic distribution of the zeros of these multiple orthogonal polynomials by using a result on the asymptotic
behavior of the ratio of two neighboring polynomials \cite{WVA-RA}. This result uses the nearest neighbor recurrence relations for multiple orthogonal
polynomials
\[   xP_{\vec{n}}(x) = P_{\vec{n}+\vec{e}_k}(x) + b_{\vec{n},k} P_{\vec{n}}(x) + \sum_{j=1}^r a_{\vec{n},j} P_{\vec{n}-\vec{e}_j}(x), \qquad
    1 \leq k \leq r, \]
where $\vec{e}_j=(0,\ldots,0,1,0,\ldots,0)$ with $1$ in the $j$th entry,
and some knowledge about the asymptotic behavior for the recurrence coefficients $a_{\vec{n},j},b_{\vec{n},j}$ $(1 \leq j \leq r)$. The ratio
asymptotic behavior for Jacobi-Pi\~neiro polynomials will be obtained in Section \ref{sec:JPratio} and for multiple Laguerre polynomials of the first kind
in Section \ref{sec:MLratio}. The asymptotic distribution of the zeros of Jacobi-Pi\~neiro polynomials will be obtained in Section \ref{sec:JPzero},
where the following result will be proved. We will use the multi-index $\vec{1}=(1,1,\ldots,1)$ so that the diagonal index is $(n,n,\ldots,n) = n\vec{1}$.

\begin{theorem}  \label{thm:JP}
Let $0 < x_{1,rn} < x_{2,rn} < \cdots < x_{rn,rn} < 1$ be the zeros of the Jacobi-Pi\~neiro polynomial $P_{n\vec{1}}$ with multi-index
$n\vec{1}=(n, n ,\ldots,  n )$. Then for every continuous function $f$ on $[0,1]$ one has
\[   \lim_{n \to \infty} \frac{1}{rn} \sum_{k=1}^{rn} f(x_{k,rn}) = \int_0^1 f(t) v_r(t)\, dt, \]
where the density $v_r$ on $[0,1]$ is given by means of a density $w_r$ on $[0,c_r]$ as
\[    v_r(x) = c_r w_r(c_rx), \qquad c_r = \frac{(r+1)^{r+1}}{r^r}, \]
and with the change of variables
\[     \hat{x} = c_r x = \frac{\bigl(\sin (r+1) \varphi \bigr)^{r+1}}{\sin \varphi \bigl( \sin r\varphi \bigr)^r}, \qquad 0 < \varphi < \frac{\pi}{r+1}, \]
the density $w_r$ is
\begin{eqnarray}   \label{wr}
   w_r(\hat{x}) &=& \frac{r+1}{\pi} \frac{1}{|\hat{x}'(\varphi)|} \nonumber \\
      &=& \frac{r+1}{\pi \hat{x}} \frac{\sin\varphi \sin r\varphi \sin(r+1)\varphi}{(r+1)^2 \sin^2 r\varphi -2r(r+1)\sin(r+1)\varphi \sin r\varphi \cos \varphi + r^2 \sin^2(r+1)\varphi}. \qquad
\end{eqnarray}
\end{theorem}

\begin{figure}[ht]
\centering
\includegraphics[width=4.5in]{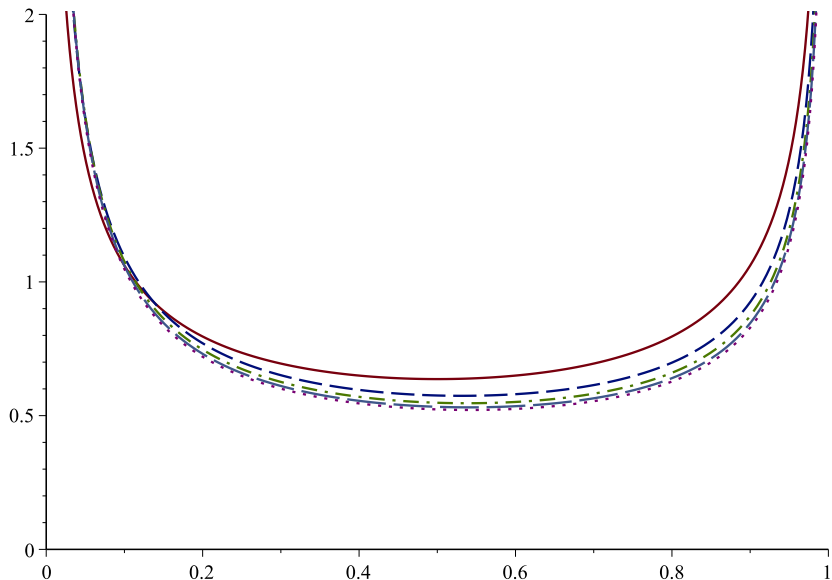}
\caption{The asymptotic zero densities $v_r$ for Jacobi-Pi\~neiro polynomials: $r=1$ (solid), 
$r=2$ (dash), $r=3$ (dash dot), $r=4$ (long dash), and $r=5$ (dots).}
\label{fig:JP}
\end{figure}

The density $w_r$ is in fact the uniform density on $[0, \frac{\pi}{r+1}]$ in the variable $\varphi$ since
\[   \int_0^{c_r} f(\hat{x}) w_r(\hat{x})\, d\hat{x} = \int_0^{\frac{\pi}{r+1}} f(\hat{x}(\varphi)) w_r(\hat{x}(\varphi)) |\hat{x}'(\varphi)|\, d\varphi
   = \frac{r+1}{\pi} \int_0^{\frac{\pi}{r+1}} f(\hat{x}(\varphi)) \, d\varphi.  \]
In this sense Theorem \ref{thm:JP} is the extension to multiple orthogonal polynomials of the equidistribution result for zeros of orthogonal polynomials
\cite[Thm. 12.7.2]{Szego} for the case of Jacobi-Pi\~neiro polynomials. In fact, the same asymptotic distribution of zeros will hold
for all multiple orthogonal polynomials for which the nearest neighbor recurrence coefficients behave as in 
\eqref{JP-alim}--\eqref{JP-blim}, provided the zeros of neighboring polynomials interlace.
We have plotted the density $v_r$ on $[0,1]$ for $1 \leq r \leq 5$ in Figure \ref{fig:JP}.

Observe that
\[    \hat{x}= c_r - \binom{r+1}{2} c_r \varphi^2 + \mathcal{O}(\varphi^4), \qquad \varphi \to 0, \]
and
\[    \hat{x} = \left( \frac{r+1}{\sin \frac{\pi}{r+1}}  \right)^{r+1} \left( \frac{\pi}{r+1} - \varphi \right)^{r+1} + \mathcal{O}\left( \left( \frac{\pi}{r+1} - \varphi \right)^{r+2} \right) , \qquad \varphi \to \frac{\pi}{r+1}, \]
so that the density $v_r$ behaves as $(\hat{x}=c_r x)$
\[  v_r(x) = \mathcal{O}(\varphi^{-1}) = \mathcal{O}\left((1-x)^{-1/2}\right), \qquad  x \to 1, \]
and
\[  v_r(x) = \mathcal{O}\left( \left( \frac{\pi}{r+1} - \varphi \right)^{-r} \right) = \mathcal{O}\left(x^{-\frac{r}{r+1}} \right), \qquad x \to 0. \]
Hence the densities $v_r$ have a square root singularity at $1$ but a higher order singularity at $0$ when $r > 1$, which means that
the zeros are more dense near the endpoints $0$ and $1$, and even more so near $0$ than near $1$ when $r > 1$. For $r=1$ the density $v_1$ is the well-known arcsin density on $[0,1]$, 
\[    v_1(x) = \frac{1}{\pi} \frac{1}{\sqrt{x(1-x)}}, \qquad 0 < x < 1, \]
which is the equilibrium measure for $[0,1]$ in logarithmic potential theory. For $r=2$ the density can explicitly be written as
\[  v_2(x) = \frac{\sqrt{3}}{4\pi} \frac{(1+\sqrt{1-x})^{1/3} + (1-\sqrt{1-x})^{1/3}}{x^{2/3} \sqrt{1-x}}, \qquad x \in (0,1), \]
and this asymptotic zero distribution was already found in \cite[Thm. 2.5]{Cous2-WVA}. The moments of $w_r$ are integers given by
\[   \int_0^{c_r} x^n w_r(x)\, dx = \frac{r+1}{\pi} \int_0^{\frac{\pi}{r+1}} x(\varphi)^n \, d\varphi =  \binom{(r+1)n}{n}, 
\qquad n \in \mathbb{N} = \{0,1,2,\ldots\},  \]
which follows from \cite[Remark 3.4]{Neuschel}.

For multiple Laguerre polynomials we need to use a scaling to prevent the zeros from going to infinity. The appropriate scaling is to divide all the
zeros of $L_{\vec{n}}$ by $|\vec{n}|$, so that we are in fact investigating the zeros of $L_{n\vec{1}}(rnx)$ for the multi-index 
$n\vec{1} = (n , n ,\ldots, n)$
on the diagonal. In Section \ref{sec:MLzero} we obtain  the asymptotic distribution of the scaled zeros, where we prove the following result. 

\begin{theorem}  \label{thm:ML}
Let $0 < x_{1,rn} < x_{2,rn} < \cdots < x_{rn,rn}$ be the zeros of the multiple Laguerre polynomials $L_{n\vec{1}}$ with multi-index
$n\vec{1}=(n,n, \ldots, n)$. Then for every continuous function $f$ on $[0,c_r/r]$ one has
\[   \lim_{n \to \infty} \frac{1}{rn} \sum_{k=1}^{rn} f\left( \frac{x_{k,rn}}{rn} \right) = \int_0^{c_r} f(t/r) u_r(t)\, dt,
\qquad c_r = \frac{(r+1)^{r+1}}{r^r},  \]
where the density $u_r$ on $[0,c_r]$ is given by
\begin{equation}  \label{ur}
    u_r(\hat{x}) = \frac{1}{r\pi} \frac{(\sin r\varphi)^{r+1}}{\bigl(\sin (r+1)\varphi \bigr)^r}, \qquad 0 < \varphi < \frac{\pi}{r+1},  
\end{equation}
where
\[     \hat{x} = \frac{\bigl(\sin (r+1) \varphi \bigr)^{r+1}}{\sin \varphi \bigl( \sin r\varphi \bigr)^r}, \qquad 0 < \varphi < \frac{\pi}{r+1}. \]
\end{theorem}

\begin{figure}[ht]
\centering
\includegraphics[width=4.5in]{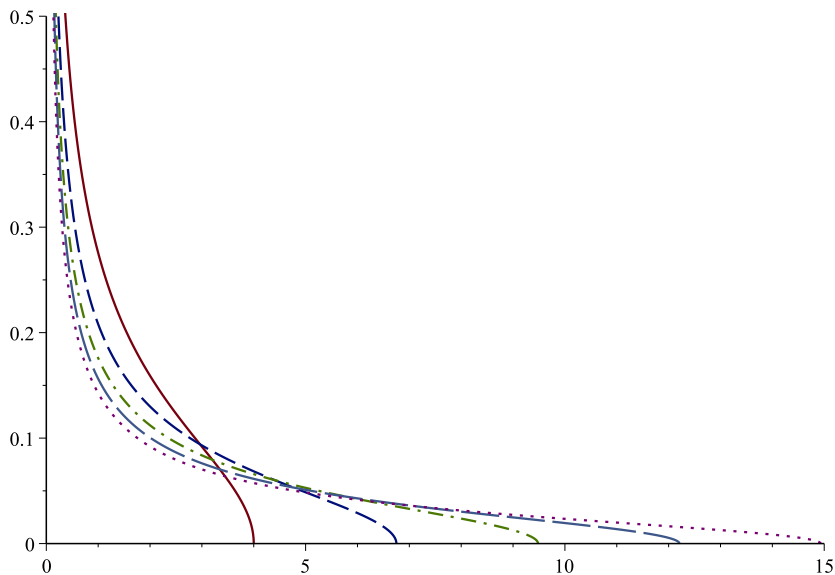}
\caption{The asymptotic zero densities $u_r$ for multiple Laguerre polynomials of the first kind: $r=1$ (solid), 
$r=2$ (dash), $r=3$ (dash dot), $r=4$ (long dash), and $r=5$ (dots).}
\label{fig:ML}
\end{figure}

The densities $u_r$ for $1 \leq r \leq 5$ are plotted in Figure \ref{fig:ML}. 

The density of the scaled zeros $\{ \frac{x_{k,rn}}{rn}, 1 \leq k \leq rn\}$
is therefore given by $ru_r(rx)$ for $0 < x < c_r/r$. Note that the densities $u_r$ behave as
\[   u_r(\hat{x}) = \mathcal{O}(\varphi) = \mathcal{O}\bigl(  (c_r - \hat{x})^{1/2} \bigr), \qquad \hat{x} \to c_r, \]
and 
\[  u_r(\hat{x}) = \mathcal{O}\left( \left(\frac{\pi}{r+1} - \varphi \right)^{-r} \right) = \mathcal{O}\left( \hat{x}^{-\frac{r}{r+1}} \right), 
\qquad \hat{x} \to 0. \]
Hence the densities $u_r$ tend to zero as a square root near the endpoint $c_r$ and have the same singularity near $0$ as in the Jacobi-Pi\~neiro case.
For $r=1$ the density is the Marchenko-Pastur density \cite{MarchenkoPastur}
\begin{equation}  \label{MP}
    u_1(\hat{x}) = \frac{1}{2\pi} \sqrt{ \frac{4-\hat{x}}{\hat{x}} }, \qquad 0 < \hat{x} < 4, 
\end{equation}
which is also the known asymptotic distribution of the (scaled) zeros of Laguerre polynomials (see, e.g., \cite{Gawronski}).
For $r=2$ we have $u_2(\hat{x}) = \frac{8}{27}g(\frac{8\hat{x}}{27})$, where
\[   g(y) = \frac{3\sqrt{3}}{16 \pi} \frac{(1+3\sqrt{1-y})(1-\sqrt{1-y})^{1/3} - (1-3\sqrt{1-y})(1+\sqrt{1-y})^{1/3}}{y^{2/3}}, \]
and the asymptotic zero distribution of the zeros of multiple Laguerre polynomials for that case was already obtained in \cite[Thm. 2.6]{Cous2-WVA}.
 An interesting observation is that the moments of $u_r$ are given by
\[  \int_0^{c_r} x^n u_r(x)\, dx = \frac{1}{n+1} \binom{(r+1)n}{n}, \qquad n \in \mathbb{N}.  \]
The simple expressions for the moments of $w_r$ and $u_r$ on $[0,c_r]$ is the main reason why we prefer to express the asymptotic zero densities
in terms of densities on $[0,c_r]$, rather than on $[0,1]$ and $[0,c_r/r]$ respectively. In Section \ref{sec:FC} we will show that these densities
and the asymptotic behavior of the ratio of Jacobi-Pi\~neiro and multiple Laguerre polynomials of the first kind are related to the Fuss-Catalan
distribution with density 
\[   g_r(x) = \frac{1}{\pi} \frac{\sin^2\varphi \bigl(\sin r\varphi \bigr)^{r-1}}{\bigl(\sin (r+1)\varphi \bigr)^r}, \qquad 0 < x < c_r, \]
where
\[     x = \frac{\bigl(\sin (r+1)\varphi \bigr)^{r+1}}{\sin \varphi \bigl( \sin r\varphi \bigr)^r} , \qquad 0 < \varphi < \frac{\pi}{r+1}, \]
for which the moments are the Fuss-Catalan numbers \cite[p.~347]{Concrete}
\[    \int_0^{c_r} x^n g_r(x)\, dx = \frac{1}{rn+1} \binom{(r+1)n}{n}, \qquad n \in \mathbb{N}.   \]

\section{Ratio asymptotics for Jacobi-Pi\~neiro polynomials} \label{sec:JPratio}

The nearest neighbor recurrence relations are
\[   xP_{\vec{n}}(x) = P_{\vec{n}+\vec{e}_k}(x) + b_{\vec{n},k} P_{\vec{n}}(x) + \sum_{j=1}^r a_{\vec{n},j} P_{\vec{n}-\vec{e}_j}(x), 
\qquad 1 \leq k \leq r, \]
where the recurrence coefficients are given by
\begin{multline} \label{JP-a} 
   a_{\vec{n},j} = \frac{n_j(n_j+\alpha_j)(|\vec{n}|+\beta)}
  {(|\vec{n}|+n_j+\alpha_j+\beta+1)(|\vec{n}|+n_j+\alpha_j+\beta)(|\vec{n}|+n_j+\alpha_j+\beta-1)}  \\
 \times  \prod_{i=1}^r\frac{|\vec{n}|+\alpha_i+\beta}{|\vec{n}|+n_i+\alpha_i+\beta}
         \prod_{i=1,i\neq j}^r \frac{n_j+\alpha_j-\alpha_i}{n_j-n_i+\alpha_j-\alpha_i}, \qquad 1 \leq j \leq r,
\end{multline}
and
\begin{multline}  \label{JP-b}
  b_{\vec{n},k} = (|\vec{n}|+\beta+1) \frac{\prod_{j=1}^r (|\vec{n}|+\beta+\alpha_j+1)}{(|\vec{n}|+n_k+\beta+\alpha_k+2)\prod_{j\neq k} (|\vec{n}|+n_j+\beta+\alpha_j+1)} \\
 - (|\vec{n}|+\beta) \frac{\prod_{j=1}^r (|\vec{n}|+\beta+\alpha_j)}{\prod_{j=1}^r (|\vec{n}|+n_j+\beta+\alpha_j)}, \qquad 1 \leq k \leq r.  
\end{multline}
(see, e.g., \cite{WVA-NN}).

If we take the multi-index $\vec{n}=(\lfloor q_1 n\rfloor , \ldots, \lfloor q_r n \rfloor)$, where $q_j > 0$, $\sum_{j=1}^r q_j = 1$, 
and $\lfloor \cdot \rfloor$ is the floor function (i.e., $\lfloor a \rfloor = k$ whenever $k \leq a < k+1$), then
the asymptotic behavior of the recurrence coefficients is
\begin{equation}  \label{JP-alim}
   \lim_{n \to \infty} a_{\vec{n},j} = \frac{q_j^{r+1}}{(1+q_j)^3} \prod_{k=1}^r \frac{1}{1+q_k} \prod_{i \neq j} \frac{1}{q_j-q_i} =: a_j, 
\qquad 1 \leq j \leq r, 
\end{equation}
and with a bit of elementary calculus
\begin{equation}  \label{JP-blim}
  \lim_{n \to \infty} b_{\vec{n},j} = \prod_{k=1}^r \frac{1}{1+q_k} \left( r+1 - \sum_{k=1}^r \frac{1}{1+q_k} - \frac{1}{1+q_j} \right) =: b_j,
  \qquad 1 \leq j \leq r.
\end{equation}
In order to have  finite values of $a_j$, we assume for the moment that $q_i \neq q_j$ whenever $i \neq j$, but later on we will take the limit
$q_j \to 1/r$ for every $j$. This passage to the limit is allowed since the asymptotic distribution of the zeros is continuous in the parameters
$(q_1,\ldots,q_r)$, which can be shown as in \cite[Thm.~2]{DomWVA}. 
We will use the notation
\[      p(\vec{q}) = \prod_{k=1}^r \frac{1}{1+q_k}, \quad s = r+1 - \sum_{k=1}^r \frac{1}{1+q_k} , \]
so that
\begin{equation}  \label{JP-ab}
   a_j = p(\vec{q}) \frac{q_j^{r+1}}{(1+q_j)^3} \prod_{i \neq j} \frac{1}{q_j-q_i}, \quad  b_j = p(\vec{q}) \left(s- \frac{1}{1+q_j} \right). 
\end{equation}

According to \cite[Thm. 1.1]{WVA-RA}, the ratio asymptotics for the Jacobi-Pi\~neiro polynomials with multi-index
$\vec{n}=(\lfloor q_1 n\rfloor , \ldots, \lfloor q_r n \rfloor)$ will then be given by
\begin{equation}   \label{JP-ratio}
   \lim_{n \to \infty} \frac{P_{\vec{n}+\vec{e}_k}(x)}{P_{\vec{n}}(x)} = z(x) - b_k, \qquad 1 \leq k \leq r, 
\end{equation}
uniformly on compact subsets of $\mathbb{C} \setminus [0,1]$, where $z$ is the solution of the algebraic equation 
\begin{equation}  \label{algebeq}
   (z-x) B_r(z) + A_{r-1}(z) = 0
\end{equation}for which $z(x) -x \to 0$ when $x \to \infty$. 
In \cite{WVA-RA} the convergence was given uniformly on compact subsets of $\mathbb{C} \setminus \mathbb{R}$, but since all the zeros
of Jacobi-Pi\~neiro polynomials are in $[0,1]$, the Stieltjes-Vitali theorem can be used to extend the uniform convergence to compact subsets
of $\mathbb{C} \setminus [0,1]$.
Here $B_r(z) = \prod_{j=1}^r (z-b_j)$ and $A_{r-1}$ is the polynomial
of degree $r-1$ for which
\[   \frac{A_{r-1}(z)}{B_r(z)} = \sum_{j=1}^r \frac{a_j}{z-b_j} .  \]
The latter means that the residue of $A_{r-1}/B_r$ at $b_j$ is given by $a_j$:
\begin{equation}    \label{JP-res}  
   a_j = \frac{A_{r-1}(b_j)}{B_r'(b_j)} = \frac{A_{r-1}(b_j)}{\prod_{i \neq j}( b_j-b_i )} . 
\end{equation}
Observe that
\[   \prod_{i \neq j} (b_j-b_i) = \frac{p(\vec{q})^r}{(1+q_j)^{r-2}} \prod_{i\neq j} (q_j-q_i), \]
so that the condition on the residues \eqref{JP-res} becomes
\begin{equation}  \label{A-interp}
   A_{r-1}(b_j) = \left( \frac{p(\vec{q}) q_j}{1+q_j} \right)^{r+1}, \qquad 1 \leq j \leq r. 
\end{equation}
This is a Lagrange interpolation problem.
If we use \eqref{JP-ab} to write $q_j$ in terms of $b_j$, then
\[   q_j = \frac{1}{s-b_j/p(\vec{q})}-1, \]
so that
\[    \frac{p(\vec{q})q_j}{1+q_j} = p(\vec{q})(1-s)+b_j, \qquad 1 \leq j \leq r. \]
The interpolation problem \eqref{A-interp} then becomes
\[  A_{r-1}(b_j) = \bigl( b_j + p(\vec{q})(1-s) \bigr)^{r+1}, \qquad 1 \leq j \leq r, \]
hence $A_{r-1}(z)$ is a polynomial of degree $r-1$ interpolating the polynomial $\bigl( z+p(\vec{q})(1-s)\bigr)^{r+1}$ at the points $b_j$ $(1 \leq j \leq r)$.
If we take the limit where $q_j \to 1/r$ for every $j$, then
\[    p(\vec{q}) \to \left( \frac{r}{r+1} \right)^r =:p, \quad s \to \frac{2r+1}{r+1}, \quad b_j \to p \left( s- \frac{r}{r+1} \right) = p, \]
hence all the interpolation points coincide. It is well known that the Lagrange interpolating polynomial for which all the interpolation points coincide
corresponds to the Taylor polynomial of degree $r-1$ of the function $f(z) = \bigl( z + p(1-s) \bigr)^{r+1}$ around the common interpolation point $p$.
This Taylor polynomial of degree $r-1$ is the polynomial $\bigl( z + p(1-s) \bigr)^{r+1}$ of degree $r+1$ from which we subtract the
last two terms of the Taylor expansion around $p$:
\begin{eqnarray}  \label{JP-Afin} 
       A_{r-1}(z) &=& \bigl( z + p(1-s) \bigr)^{r+1} - (z-p)^{r+1} \frac{f^{(r+1)}(p)}{(r+1)!} - (z-p)^r \frac{f^{(r)}(p)}{r!} \nonumber \\
                  &=&  \bigl( z + p(1-s) \bigr)^{r+1} - (z-p)^{r+1} - (r+1) p(2-s) (z-p)^r .    
\end{eqnarray}
The algebraic equation \eqref{algebeq} for multi-indices on the diagonal then becomes
\[     (z-x) (z-p)^r + \bigl( z + p(1-s) \bigr)^{r+1} - (z-p)^{r+1} - (r+1)p(2-s)(z-p)^r = 0, \]
which simplifies to
\begin{equation}  \label{JP-algeqfin}
     x (z-p)^r = \Bigl( z - \frac{pr}{r+1}  \Bigr)^{r+1} .
\end{equation}

\section{Relation with the Fuss-Catalan numbers}  \label{sec:FC}

Recently the Fuss-Catalan distribution and other related distributions (Raney distributions) appeared as limiting distributions
of eigenvalues and singular values of certain random matrices \cite{ForLiu}, \cite{PenZycz}, \cite{Neuschel}. In this section we will show
how the ratio asymptotics in \eqref{JP-ratio} is related to the Stieltjes transform of the Fuss-Catalan distribution.
The weights $w_r$ and $u_r$ in Theorem \ref{thm:JP} and \ref{thm:ML} cannot be identified with the Fuss-Catalan distribution or any of the Raney
distributions (except $u_1$, which is the Catalan distribution) because their behavior near the endpoints of the interval differs from
the behavior of the Raney distributions given in \cite{MlotPenZycz}.

If we scale the variables $\hat{x}=c_r x$ and $\hat{z}=c_r z$, where
\[    c_r = \frac{(r+1)^{r+1}}{r^r} = \frac{r+1}{p}, \]
then the algebraic equation \eqref{JP-algeqfin} becomes
\begin{equation}  \label{JP-algeqhat}
  \hat{x} (\hat{z} - r-1 )^r = ( \hat{z} - r )^{r+1}.  
\end{equation}
If we define 
\begin{equation}  \label{w-z}
    \omega = \frac{\hat{z}-r}{\hat{z}-r-1},  \quad \hat{z} = \frac{(r+1)\omega-r}{\omega-1}, 
\end{equation}
then the algebraic equation becomes  
\begin{equation}  \label{algeqFC}
   \omega^{r+1} + \hat{x} - \hat{x} \omega = 0.  
\end{equation}
This is the algebraic equation for the generating function $G(1/\hat{x})$ of the Fuss-Catalan numbers \cite[p.~347]{Concrete} \cite[Eq. (3.12)]{Neuschel}.
As in \cite[\S 3]{Neuschel}, we assume that a solution exists of the form $\omega = \rho e^{i\varphi}$, where $\rho >0$ and $\varphi$ is real. Then
inserting this in \eqref{algeqFC} gives
\[   \rho^{r+1} e^{i(r+1)\varphi} + \hat{x} - \hat{x} \rho e^{i\varphi} = 0. \]
This gives for the real and the imaginary part
\begin{eqnarray}
      \rho^{r+1} \cos(r+1)\varphi + \hat{x} - \hat{x} \rho \cos \varphi &=& 0,   \label{rho1} \\ 
      \rho^{r+1} \sin(r+1)\varphi - \hat{x} \rho \sin \varphi  &=& 0.   \label{rho2} 
\end{eqnarray}
From \eqref{rho2} we find
\begin{equation}   \label{x-rho}
      \hat{x} = \rho^r \frac{\sin (r+1)\varphi}{\sin \varphi}, 
\end{equation}
and inserting this in \eqref{rho1} gives
\begin{equation}   \label{rho}
      \rho(\hat{x}) = \frac{\sin(r+1)\varphi}{\sin r \varphi},
\end{equation}
from which
\begin{equation}  \label{x-phi}
    \hat{x} = \frac{\bigl( \sin (r+1)\varphi \bigr)^{r+1}}{\sin \varphi \bigl( \sin r\varphi \bigr)^{r}}.
\end{equation}
Observe that $\rho(x) >0$ for $0 < \varphi < \frac{\pi}{r+1}$, and $\hat{x}$ is a monotonically decreasing function mapping $[0, \frac{\pi}{r+1}]$
into $[0,c_r]$. So for $\hat{x} \in [0,c_r]$ there is a solution $\rho e^{i\varphi}$ of the algebraic equation \eqref{algeqFC}. The conjugate 
function $\rho e^{-i\varphi}$ is also a solution for $\hat{x} \in [0,c_r]$. In fact both solutions are the boundary value of a function
$\omega$ which is analytic on $\mathbb{C} \setminus [0,c_r]$ and
\begin{equation}  \label{wpm}
     \omega_+ = \lim_{\epsilon \to 0+} \omega(\hat{x}+i\epsilon) = \rho e^{-i\phi}, \quad   
     \omega_- = \lim_{\epsilon \to 0+} \omega(\hat{x}-i\epsilon) = \rho e^{i\phi}, 
\end{equation}
because this $\omega$ is $G(1/\hat{x}) = \hat{x} F(\hat{x})$, where $F$ is the Stieltjes transform of the Fuss-Catalan distribution 
\[    F(z) = \int_0^{c_r} \frac{g_r(y)}{z-y}\, dy, \qquad z \in \mathbb{C} \setminus [0,c_r], \]
with $g_r$ the Fuss-Catalan density, and a Stieltjes transform has the property that
\[     \textup{Im\,} F(z) \begin{cases}  < 0, & \textup{Im\,} z >0, \\
                                         > 0, & \textup{Im\,} z < 0.   \end{cases}.  \]
 
Observe that
\[   \frac{1}{z-p} = \lim_{n \to \infty} \frac{P_{\vec{n}}(x)}{P_{\vec{n}+\vec{e}_k}(x)} \]
is the Stieltjes transform of a probability measure on $[0,1]$, since we have 
\[   \frac{P_{\vec{n}}(x)}{P_{\vec{n}+\vec{e}_k}(x)} = \sum_{j=1}^{|\vec{n}|+1} \frac{c_{j,\vec{n}}}{x-x_{j,\vec{n}+\vec{e}_k}}, \]
and $c_{j,\vec{n}} > 0$ because the zeros of $P_{\vec{n}}$ and $P_{\vec{n}+\vec{e}_k}$ interlace \cite[Thm.2.1]{HanWVA}, and
$\sum c_{j,\vec{n}} = 1$ since we are dealing with monic polynomials. With the change of variables
$\hat{x}=c_rx$ and $\hat{z} = c_rz$ it follows that 
\[ 1/(\hat{z}-r-1) = \int_0^{c_r} \frac{d\mu(y)}{\hat{x}-y}  \]
is the Stieltjes transform of a probability distribution $\mu$ on $[0,c_r]$. Note that \eqref{w-z} implies
\[   \frac{1}{\hat{z}-r-1} = \omega-1 = \hat{x}F(\hat{x})-1, \]
where $F$ is the Stieltjes transform of the Fuss-Catalan distribution,
\[   F(\hat{x}) = \frac{1}{\hat{x}} \sum_{n=0}^\infty  \frac{1}{rn+1} \binom{(r+1)n}{n} \frac{1}{\hat{x}^n}, \]
so that $1/(\hat{z}-r-1)$ is the Stieltjes transform of the probability measure for which the moments are the Fuss-Catalan numbers shifted by one
\[   \int_0^{c_r} y^n \, d\mu(y) = \frac{1}{r(n+1)+1} \binom{(r+1)(n+1)}{n+1},  \]
and hence this probability distribution has a density $\hat{x} g_r(\hat{x})$, where $g_r$ is the Fuss-Catalan density
\[   g_r(\hat{x}) = \frac{1}{\pi} \frac{\sin^2 \varphi ( \sin r \varphi )^{r-1}}{\bigl( \sin (r+1)\varphi \bigr)^r},  
  \qquad 0 < \varphi < \frac{\pi}{r+1}, \]
where $\hat{x}$ is given in \eqref{x-phi}. In particular this gives
\[    \frac{1}{\hat{z} - r -1} = \int_0^{c_r} \frac{yg_r(y)}{\hat{x}-y}  \, dy.  \]
The weight is explicitly given by
\[   \hat{x}g_r(\hat{x}) = \frac{1}{\pi} \frac{\sin \varphi \sin (r+1)\varphi}{\sin r \varphi}, \qquad 0 \leq \varphi < \frac{\pi}{r+1}, \]
with $\hat{x}$ as in \eqref{x-phi}.

\section{Proof of Theorem \ref{thm:JP}}  \label{sec:JPzero}

So far we found that for $\vec{n}$ near the diagonal (i.e., $n_j/n \to 1/r$ for every $j$) one has
\begin{equation}  \label{ratio}
  \lim_{n \to \infty} \frac{P_{\vec{n}+\vec{e}_k}(x)}{P_{\vec{n}}(x)} = z(x) - p =  \frac{1}{c_r} (\hat{z} - r- 1), 
\end{equation}
uniformly for $x$ on compact subsets of $\mathbb{C} \setminus [0,1]$, or $\hat{x}$ on compact subsets of
$\mathbb{C} \setminus [0,c_r]$. However we are interested in the asymptotic behavior of
\[   \frac{1}{|\vec{n}|} \frac{P_{\vec{n}}'(x)}{P_{\vec{n}}(x)}, \]
where the prime $'$ denotes the derivative with respect to $x$,
because the limit will give the Stieltjes transform of the asymptotic distribution of the zeros of $P_{\vec{n}}$.
By taking derivatives with respect to $x$ in \eqref{ratio} we find
\[   \lim_{n \to \infty} \frac{P_{\vec{n}+\vec{e}_k}(x)}{P_{\vec{n}}(x)}
   \left(  \frac{P_{\vec{n}+\vec{e}_k}'(x)}{P_{\vec{n}+\vec{e}_k}(x)} - \frac{P_{\vec{n}}'(x)}{P_{\vec{n}}(x)} \right)  = z' = \frac{\hat{z}'}{c_r}, \]
uniformly for $x$ on compact subsets of $\mathbb{C} \setminus [0,1]$. Together with \eqref{ratio}, this gives
\[  \lim_{n \to \infty} \left( \frac{P_{\vec{n}+\vec{e}_k}'(x)}{P_{\vec{n}+\vec{e}_k}(x)} - \frac{P_{\vec{n}}'(x)}{P_{\vec{n}}(x)} \right)
      = \frac{\hat{z}'}{\hat{z}-r-1}. \]
If we use this result successively for each $k$, $1 \leq k \leq r$, then we find for multi-indices on the diagonal
$n\vec{1} = (n,n,\ldots,n)$ and $(n+1)\vec{1} = (n+1,n+1,\ldots,n+1)$
\[   \lim_{n \to \infty} \left( \frac{P_{(n+1)\vec{1}}'(x)}{P_{(n+1)\vec{1}}(x)} - \frac{P_{n\vec{1}}'(x)}{P_{n\vec{1}}(x)} \right)
    = \frac{r\hat{z}'}{\hat{z}-r-1}. \] 
Then by taking averages (Ces\`aro's lemma) we get
\[   \lim_{n \to \infty} \frac{1}{n} \sum_{k=0}^{n-1} \left( \frac{P_{(k+1)\vec{1}}'(x)}{P_{(k+1)\vec{1}}(x)} 
 - \frac{P_{k\vec{1}}'(x)}{P_{k\vec{1}}(x)} \right) = \frac{r\hat{z}'}{\hat{z}-r-1}, \]
and since this contains a telescoping sum, this becomes
\[  \lim_{n  \to \infty}  \frac{1}{rn} 
 \frac{P_{n\vec{1}}'(x)}{P_{n\vec{1}}(x)} = \frac{\hat{z}'}{\hat{z}-r-1}, \]
uniformly for $x$ on compact subsets of $\mathbb{C} \setminus [0,1]$, so that the right hand side is the Stieltjes transform
of the asymptotic zero distribution of the zeros of $P_{n\vec{1}}$ for $n\vec{1} = (n,n,\ldots,n)$.
From the Stieltjes transform we can find the density by using Stieltjes' inversion formula
\[    2\pi i v_r(x) =  \left(\frac{\hat{z}'}{\hat{z}-r-1}\right)_- - \left(\frac{\hat{z}'}{\hat{z}-r-1}\right)_+ . \]
Taking derivatives in \eqref{JP-algeqhat} (and recalling that $\hat{x}=c_r x$) gives
\[     c_r (\hat{z}-r-1)^r + \hat{x} r(\hat{z}- r- 1)^{r-1} \hat{z}' = (r+1) \left( \hat{z}- r \right)^r \hat{z}', \]
so that
\[   \frac{\hat{z}'}{\hat{z}-r-1} = \frac{c_r}{\hat{x}} \frac{\hat{z}- r}{\hat{z}-2r-1} .  \]
Writing this in terms of $\omega$ using \eqref{w-z} gives
\[   \frac{\hat{z}'}{\hat{z}-r-1} = \frac{c_r}{\hat{x}} \frac{\omega}{-r\omega+r+1}. \]
Now use $\omega_+ = \rho e^{-i\varphi}$ and $\omega_- = \rho e^{i\varphi}$ to find the density
\[     v_r(x) = \frac{r+1}{\pi x} \frac{\rho \sin \varphi}{|r\rho e^{i\varphi} -r-1|^2},  \]
and clearly $v_r(x)=v_r(\hat{x}/c_r) = c_r w_r(\hat{x})$, with the weight in \eqref{wr}. 
Observe that $\hat{x}: [0,\frac{\pi}{r+1}] \to [0,c_r]$ is a monotonically decreasing function with
\[  \hat{x}'(\varphi) = \frac{-\hat{x}}{\sin \varphi \sin r\varphi \sin (r+1) \varphi} |(r+1)\sin r\varphi - e^{i\varphi} r \sin(r+1)\varphi|^2 \]
so that
\[  w_r(\hat{x}) =  \frac{r+1}{\pi} \frac{1}{|\hat{x}'(\varphi)|}, \qquad 0 < \varphi < \frac{\pi}{r+1}. \]

\section{Ratio asymptotics for multiple Laguerre polynomials of the first kind}  \label{sec:MLratio}

The nearest neighbor recurrence relations for multiple orthogonal polynomials of the first kind are given by
\[   xL_{\vec{n}}(x) = L_{\vec{n}+\vec{e}_k}(x) + b_{\vec{n},k} L_{\vec{n}}(x) + \sum_{j=1}^r a_{\vec{n},j} L_{\vec{n}-\vec{e}_j}(x), 
\qquad 1 \leq k \leq r, \]
where the recurrence coefficients are given by
\begin{equation}   \label{ML-a}
  a_{\vec{n},j} = n_j(n_j+\alpha_j) \prod_{i=1,i\neq j}^r \frac{n_j+\alpha_j-\alpha_i}{n_j-n_i+\alpha_j-\alpha_i}, \qquad 1 \leq j \leq r, 
\end{equation}
and
\begin{equation}  \label{ML_b}
   b_{\vec{n},k} =  |\vec{n}|+n_k+\alpha_k+1, \qquad 1 \leq k \leq r.   
\end{equation}  
(see, e.g., \cite{WVA-NN}).
We can now proceed as in the case of Jacobi-Pi\~neiro polynomials. The recurrence coefficients are somewhat easier but they are unbounded
so that we need to use a scaling. Suppose again that $\vec{n} = (\lfloor q_1n \rfloor, \ldots, \lfloor q_r n \rfloor )$, where $q_i\neq q_j$
whenever $i \neq j$. It then follows that
\[  \lim_{n \to \infty}  \frac{a_{\vec{n},j}}{n^2} = q_j^{r+1} \prod_{i \neq j} \frac{1}{q_j-q_i} =: a_j, \qquad 1 \leq j \leq r, \]
and 
\[  \lim_{n \to \infty}  \frac{b_{\vec{n},j}}{n} = 1+q_j =: b_j, \qquad 1 \leq j \leq r. \]
According to \cite[Thm. 1.2]{WVA-RA} we then have
\begin{equation}  \label{ML-ratio}
   \lim_{n \to \infty} \frac{L_{\vec{n}+\vec{e}_k}(nx)}{nL_{\vec{n}}(nx)} = z(x) - b_k, \qquad 1 \leq k \leq r, 
\end{equation}
uniformly on compact subsets of $\mathbb{C} \setminus [0,\infty)$, where $z$ is the solution of the algebraic equation
\[    (z-x) B_r(z) + A_{r-1}(z) = 0, \]
where $B_r(z) = \prod_{j=1}^r (z-b_j)$ and $A_{r-1}$ is obtained from
\[   \frac{A_{r-1}(z)}{B_r(z)} = \sum_{j=1}^r \frac{a_j}{z-b_j} . \]
The uniform convergence on compact subsets of $\mathbb{C} \setminus \mathbb{R}$ in \cite{WVA-RA} can be extended to $\mathbb{C} \setminus [0,\infty)$
because the zeros of multiple Laguerre polynomials of the first kind are on $[0,\infty)$. One can even extend this further to $\mathbb{C} \setminus
[0,c_r/r]$ since all the scaled zeros are dense on $[0,c_r/r]$, but we will not need this here.
Observe that
\[   \prod_{i \neq j} (q_j-q_i) = \prod_{i \neq j} (b_j-b_i), \]
so that we get the interpolation condition
\[   A_{r-1}(b_j) = q_j^{r+1} = (b_j - 1)^{r+1}, \qquad 1 \leq j \leq r.  \]
Hence $A_{r-1}(z)$ is the Lagrange interpolating polynomial of degree $r-1$ for the function $f(z) = (z-1)^{r+1}$ for the interpolation 
points $b_1,\ldots, b_r$. Now let $q_j \to 1/r$ for every $j$, then 
\[  b_j \to \frac{r+1}{r}, \]
and $A_{r-1}(z)$ will be the Taylor polynomial of degree $r-1$ around $\frac{r+1}{r}$ for the function $f(z) = (z-1)^{r+1}$. This gives
\[   A_{r-1}(z) = (z - 1)^{r+1} - \left( z- \frac{r+1}{r} \right)^{r+1} - \frac{r+1}{r} \left( z - \frac{r+1}{r} \right)^{r} . \]
The algebraic equation for multi-indices near the diagonal then becomes
\begin{equation}  \label{ML-algeqfin}
     x \left( z - \frac{r+1}{r} \right)^r =  (z-1)^{r+1}.  
\end{equation}
The change of variables $rz = \hat{z}$ and $rx =\hat{x}$ gives the same algebraic equation as in \eqref{JP-algeqhat}.

\section{Proof of Theorem \ref{thm:ML}}  \label{sec:MLzero}

As in Section \ref{sec:JPzero} we use 
\[   \frac{L_{n\vec{1}}'(x)}{L_{n\vec{1}}(x)} = \sum_{k=0}^{n-1} \left( \frac{L_{(k+1)\vec{1}}'(x)}{L_{(k+1)\vec{1}}(x)}  
 - \frac{L_{k\vec{1}}'(x)}{L_{k\vec{1}}(x)}
   \right), \]
where $k\vec{1} = (k,k, \ldots, k)$ and $(k+1)\vec{1} = (k+1,k+1,\ldots,k+1)$.
However, because of the scaling, we need to consider (observe that $|n\vec{1}| = rn$)
\[   \frac{L_{n\vec{1}}'(rnx)}{rnL_{n\vec{1}}(rnx)} = \sum_{k=0}^{n-1} \left( \frac{L_{(k+1)\vec{1}}'(rnx)}{rnL_{(k+1)\vec{1}}(rnx)} 
 - \frac{L_{k\vec{1}}'(rnx)}{rnL_{k\vec{1}}(rnx)}  \right), \]
so that we can not use Ces\`aro's lemma to get the asymptotic behavior. We modify the proof as follows. For $\frac{k}{n} \leq t < \frac{k+1}{n}$ one has
$\lfloor nt \rfloor = k$, hence the sum can be written as an integral
\[  \frac{L_{n\vec{1}}'(rnx)}{rnL_{n\vec{1}}(rnx)}  
 = n \int_0^1 \left( \frac{L_{(\lfloor nt \rfloor+1) \vec{1}}'(rnx)}{rnL_{(\lfloor nt \rfloor +1)\vec{1}}(rnx)} 
 - \frac{L_{\lfloor nt \rfloor \vec{1}}'(rnx)}{rnL_{\lfloor nt \rfloor \vec{1}}(rnx)}  \right)\, dt,  \]
and the integrand can be written as
\[   \left( \frac{L_{(\lfloor nt \rfloor+1)\vec{1}}(rnx)}{rnL_{\lfloor nt \rfloor \vec{1}}(rnx)} \right)'/
    \left( \frac{L_{(\lfloor nt \rfloor+1)\vec{1}}(rnx)}{rnL_{\lfloor nt \rfloor \vec{1}}(rnx)} \right). \]
So we need to know the asymptotic behavior of the ratio
\[ \lim_{n \to \infty} \frac{L_{(\lfloor nt \rfloor+1)\vec{1}}(rnx)}{rnL_{\lfloor nt \rfloor \vec{1}}(rnx)} .  \]
If we change $n$ to $rn$ in Section \ref{sec:MLratio} then for $q_j\to\frac{1}{r}$ $(1 \leq j \leq r)$ we get the multi-index $n\vec{1}=(n,n,\ldots,n)$
and \eqref{ML-ratio} becomes
\[  \lim_{n \to \infty} \frac{L_{n\vec{1}+\vec{e}_k}(rnx)}{rnL_{n\vec{1}}(rnx)} = z(x) - \frac{r+1}r, \]
but we need to extend this for multi-indices containing the parameter $0 < t \leq 1$.
For this we need to
use the following asymptotic behavior of the recurrence coefficients: if $\vec{n} = (\lfloor nq_1 \rfloor, \lfloor nq_2 \rfloor, \ldots, \lfloor nq_r \rfloor)$ and $\vec{m} = (\lfloor ntq_1 \rfloor, \lfloor ntq_2 \rfloor, \ldots, \lfloor ntq_r \rfloor)$, then
\[  \lim_{n \to \infty}  \frac{a_{\vec{m},j}}{n^2} = t^2q_j^{r+1} \prod_{i \neq j} \frac{1}{q_j-q_i} = t^2a_j, 
\qquad 1 \leq j \leq r, \]
and 
\[  \lim_{n \to \infty}  \frac{b_{\vec{m},j}}{n} = t(1+q_j) = tb_j, \qquad 1 \leq j \leq r. \]
The required asymptotic behavior is then for $0 < t \leq 1$
\begin{equation}  \label{ML-ratio-t}
\lim_{n \to \infty}  \frac{L_{\vec{m}+\vec{e}_k}(nx)}{nL_{\vec{m}}(nx)} = z(x,t) - t b_k,
\end{equation}
uniformly for $x$ on compact subsets of $\mathbb{C} \setminus [0,\infty)$, where $z(x,t)$ satisfies the algebraic equation
\[      \bigl( z(x,t) -x \bigr) B_r(z,t) + A_{r-1}(z,t)  = 0, \]
with $B_r(z,t) = \prod_{j=1}^r (z-tb_j) = t^r B_r(z/t)$ and 
\[  \frac{A_{r-1}(z,t)}{B_r(z,t)} = \sum_{j=1}^r \frac{t^2 a_j}{z-tb_j} , \]
so that $A_{r-1}(z,t) = t^{r+1} A_{r-1}(z/t)$. Here we used $A_{r-1}(z) = A_{r-1}(z,1)$ and $B_r(z) = B_r(z,1)$, which are the polynomials
in Section \ref{sec:MLratio}.
If $q_j \to \frac{1}{r}$ $(1 \leq j \leq r)$ then $b_j \to \frac{r+1}{r}$ $(1 \leq j \leq r)$ and the algebraic equation
for $z(x,t)$ becomes
\begin{equation}   \label{algeqMLt}
  x  \left( z(x,t) - t\frac{r+1}{r} \right)^r = \bigl( z(x,t)-t \bigr)^{r+1}. 
\end{equation}
Now change $n$ to $rn$ so that we can deal with the multi-index $n\vec{1}=(n,n,\ldots,n)$.
By going from the multi-index $n\vec{1} = (n,n,\ldots,n)$ to $(n+1)\vec{1} = (n+1,n+1,\ldots,n+1)$ in $r$ steps (each time increasing one coefficient)
we then get
\[  \lim_{n \to \infty}  \frac{L_{(\lfloor nt \rfloor+1)\vec{1}}(rnx)}{(rn)^rL_{\lfloor nt \rfloor \vec{1}}(rnx)} 
= \left( z(x,t) - t \frac{r+1}{r} \right)^r , \]
so that
\[   \lim_{n \to \infty} \frac{1}{rn} \frac{L_{n\vec{1}}'(rnx)}{rnL_{n\vec{1}}(rnx)}
    = \frac{1}{r} \int_0^1 \frac{ \frac{d}{dx} \left( z(x,t) - t \frac{r+1}{r} \right)^r}{\left( z(x,t) - t \frac{r+1}{r} \right)^r}\, dt 
    = \int_0^1 \frac{z'(x,t)}{z(x,t)-t \frac{r+1}{r}}\, dt ,  \]
uniformly on compact subsets of $\mathbb{C} \setminus [0,\infty)$, where the prime $'$ means the derivative with respect to $x$.
This limit is the Stieltjes transform of the asymptotic zero distribution
\[   \int_0^{c_r/r} r\frac{u_r(rs)}{x-s} \, ds = \int_0^{c_r} \frac{u_r(y)}{x-y/r}\, dy ,    \]
and hence we have
\[    \int_0^1 \frac{z'(x,t)}{z(x,t)-t \frac{r+1}{r}}\, dt = \int_0^{c_r} \frac{u_r(y)}{x-y/r}\, dy .  \]
Observe that the change of variables $rz = t \hat{z}$ and $rx = t \hat{x}$  transforms the algebraic equation \eqref{algeqMLt}
to \eqref{JP-algeqhat}, so that $z(t\hat{x}/r,t) = t \hat{z}(\hat{x})/r$. From our analysis in Sections \ref{sec:JPratio}--\ref{sec:JPzero} we found that
\[    \frac{\hat{z}'}{\hat{z}- r-1} = \int_0^{c_r} \frac{w_r(s)}{\hat{x}-s}\, ds, \]
hence
\[    \frac{z'(x,t)}{z(x,t)-t \frac{r+1}{r}} = \frac{r}{t} \int_0^{c_r} \frac{w_r(s)}{\frac{rx}{t}-s}\, ds 
  = \int_0^{c_r} \frac{w_r(s)}{x-\frac{ts}{r}}\, ds. \]
Therefore
\begin{eqnarray*}
   \int_0^1 \frac{z'(x,t)}{z(x,t)-t \frac{r+1}{r}}\, dt &=& \int_0^1  \int_0^{c_r} \frac{w_r(s)}{x-\frac{ts}{r}}\, ds\, dt \\
                       &=& \int_0^1 \int_0^{tc_r} \frac{w_r(y/t)}{x-\frac{y}{r}} \frac{dy}{t}\, dt \\
                       &=& \int_0^{c_r} \frac{1}{x-\frac{y}{r}} \int_{y/c_r}^1 w_r(y/t)\, \frac{dt}{t} \, dy, 
\end{eqnarray*}
where we used the change of variables $ts = y$ in the second equality and Fubini's theorem for the third equality.
This means that 
\begin{equation}  \label{ur-wr}
   u_r(y) = \int_{y/c_r}^1 w_r(y/t)\, \frac{dt}{t} = \int_y^{c_r} w_r(x) \, \frac{dx}{x}, 
\end{equation}
and hence the asymptotic density of the scaled zeros $u_r$ is the Mellin convolution of the density $w_r$ given in \eqref{wr} and the uniform
distribution on $[0,1]$. This immediately gives the moments
\[  \int_0^{c_r} y^n u_r(y)\, dy = \int_0^{c_r} x^n w_r(x)\, dx  \int_0^1 t^n \,dt  = \frac{1}{n+1} \binom{(r+1)n}{n}.  \]
We still need to show that the density $u_r$ is given by the expression in \eqref{ur}. Observe that the derivative of \eqref{ur} with respect to 
$\varphi$ is
\[    \frac{r+1}{\pi} \frac{\sin \varphi (\sin r\varphi)^r}{\bigl( \sin (r+1)\varphi \bigr)^{r+1}} = \frac{r+1}{\pi \hat{x}}. \]
On the other hand, taking the derivative in \eqref{ur-wr} with respect to $\varphi$ gives
\[   \frac{du_r(\hat{x})}{d\varphi} = - \frac{w_r(\hat{x})}{\hat{x}} \hat{x}' = \frac{r+1}{\pi \hat{x}}, \]
where we used \eqref{wr} for the last equality. Thus, using $u_r(c_r)=0$, we find that
\[   u_r(\hat{x}) = \frac{1}{r\pi} \frac{ (\sin r\varphi )^{r+1}}{\bigl( \sin (r+1)\varphi \bigr)^r} .  \] 

\section{Concluding remarks}

There is yet another family of multiple orthogonal polynomials for which the asymptotic distribution of the zeros is of the same flavor.
These are multiple orthogonal polynomials associated with Meijer G-functions, which appear in the study of products of Ginibre random matrices
\cite{KuijlZhang}. The polynomials on the stepline for $|\vec{n}| = n$ are given by
\[  P_{n}(x) = (-1)^n \prod_{j=1}^r (n+\nu_j)! \ \sum_{k=0}^n \binom{n}{k} \frac{(-x)^k}{(k+\nu_1)! \cdots (k+\nu_r)!}, \]
and the asymptotic distribution of the scaled zeros $\{x_{k,n}/n^r, 1 \leq k \leq n\}$ is given in \cite[Thm. 3.2]{Neuschel}. The density is
\[     g_r(x) = \frac{1}{\pi} \frac{\sin^2\varphi (\sin r\varphi)^{r-1}}{\bigl( \sin (r+1) \varphi \bigr)^r}, \]
where again
\begin{equation}  \label{CRx}
   x = \frac{\bigl( \sin(r+1)\varphi \bigr)^{r+1}}{\sin \varphi (\sin r\varphi)^r}, \qquad 0 < \varphi < \frac{\pi}{r+1}. 
\end{equation}
This is the density of the Fuss-Catalan distribution, for which the moments are the Fuss-Catalan numbers
\[   \int_0^{c_r} x^n g_r(x)\, dx = \frac{1}{rn+1} \binom{(r+1)n}{n}, \qquad n \in \mathbb{N}.  \]
Observe that the density $g_r$ is a Mellin convolution of the density $w_r$ in \eqref{wr} and the beta($\frac{1}{r},1$) density:
\begin{eqnarray}  \label{Mellingr}
  g_r(y) & = & \frac{1}{r} \int_{y/c_r}^1  w_r(y/t) t^{1/r-1} \, \frac{dt}{t}  \nonumber \\ 
         & = &  \frac{1}{r} \int_y^{c_r}  w_r(x) \left( \frac{y}{x} \right)^{1/r-1} \, \frac{dx}{x}  \nonumber  \\ 
         & = & y^{1/r-1} \frac{r+1}{r\pi} \int_0^\theta \frac{d\varphi}{x^{1/r}} , 
\end{eqnarray}
where
\[    y(\theta) = \frac{\bigl( \sin(r+1)\theta \bigr)^{r+1}}{\sin \theta (\sin r\theta)^r}. \]
This can most easily be seen from
\[    \frac{d}{d\varphi} \left( \frac{ (\sin \varphi )^{1/r+1}}{\bigl( \sin (r+1) \varphi \bigr)^{1/r}} \right)
     = \frac{r+1}{r} \frac{1}{x^{1/r}} ,  \]
with $x$ given in \eqref{CRx},
which enables a straightforward computation of the last integral in \eqref{Mellingr}.
The case $r=1$ corresponds to the asymptotic zero distribution of Laguerre polynomials (the Marchenko-Pastur distribution \eqref{MP}).
The case $r=2$ was obtained earlier in \cite{Cous2-WVA} and corresponds to multiple orthogonal polynomials for modified Bessel functions
$K_\nu$ and $K_{\nu+1}$. The weight is then explicitly given by $g_2(x) = \frac{4}{27} h(\frac{4x}{27})$, where
\[    h(y) = \frac{3\sqrt{3}}{4\pi} \frac{(1+\sqrt{1-y})^{1/3} - (1-\sqrt{1-y})^{1/3}}{y^{2/3}}, \qquad 0 < y < 1.  \]

\section*{Acknowledgements}

This research was supported by KU Leuven research grant OT/12/073,  FWO research grant G.0934.13 and the Belgian Interuniversity Attraction Poles Programme P7/18. Thorsten Neuschel is a Research Associate (charg\'e de recherches) of FRS-FNRS (Belgian Fund for Scientific Research).

\begin{quote}
Walter Van Assche \\
Department of Mathematics \\
KU Leuven \\
Celestijnenlaan 200 B box 2400 \\
BE-3001 Leuven \\
Belgium \\
\texttt{walter@wis.kuleuven.be}
\end{quote}

\begin{quote}
Thorsten Neuschel \\
current address: \\
IRMP \\
Universit\'e Catholique de Louvain \\
Chemin du Cyclotron 2 \\
BE-1348 Louvain-la-Neuve \\
Belgium \\
\texttt{thorsten.neuschel@uclouvain.be}
\end{quote}

\end{document}